\documentclass[11pt]{paper}

\usepackage{amssymb,amsmath,enumerate,theorem,epsfig,rotating,graphics,changebar,eepic}%t1enc
\usepackage{amsfonts, booktabs}
\usepackage{a4,latexsym,parskip}
\usepackage{hyperref}
\hypersetup{pdfauthor=MS} \hypersetup{pdftitle=Fields}

\newcommand{\PP}{\mathbb{P}}
\newcommand{\A}{\mathbb{A}}

\newcommand{\C} [1][]{\mathbb{C}^{#1}}
\newcommand{\Q} {\mathbb{Q}}
\newcommand{\N} [1][] {\mathbb{N}_{#1}}
\newcommand{\R} {\mathbb{R}}
\newcommand{\F}{\mathbb{F}}
\newcommand{\Z}{\mathbb{Z}}

\newcommand{\OO}{\mathcal{O}}

\newcommand{\NS}{\mathop{\rm NS}\nolimits}
\newcommand{\Br}{\mathop{\rm Br}\nolimits}

\theoremstyle{break} \newtheorem{Theorem}{Theorem}

\newtheorem{Example}[Theorem]{Example}

\newtheorem{Definition}[Theorem]{Definition}

\newtheorem{Remark}[Theorem]{Remark}

\newtheorem{Conjecture}[Theorem]{Conjecture}

\newtheorem{Criterion}[Theorem]{Criterion}

\begin{document}
\setlength{\unitlength}{1cm}

\title{Arithmetic of K3 surfaces}

\author{Matthias Sch\"utt}

\maketitle

\abstract{We review recent developments in the arithmetic of K3 surfaces. Our focus lies on aspects of modularity, Picard number and rational points. Throughout we emphasise connections to geometry.}



\section{Introduction}
\label{s:intro}

K3 surfaces have established themselves as a connection between various areas of mathematics, as diverse as algebraic geometry, differential geometry, dynamics, number theory and string theory.
In this survey, we will focus on arithmetic aspects and investigate their deep interplay with geometry.

The arithmetic of curves is fairly well understood, but in higher dimensions much less is known. 
It is this context which makes us turn to K3 surfaces.
In 2004, Swinnerton-Dyer stated that "K3 surfaces are the simplest kind of variety about whose number-theoretic properties very little is known" \cite{S-D}.

Since 2004, there have been interesting developments in the arithmetic of K3 surfaces which this paper will review. Namely we will be concerned with the following topics:
\begin{enumerate}[1.]
\item Modularity and singular K3 surfaces;
\item Picard number and Galois action on the N\'eron-Severi group;
\item Rational points and potential density.
\end{enumerate}

The survey requires only very basic knowledge of algebraic geometry and number theory. The concepts involved will be introduced whenever they are needed, so that the motivation becomes clear. Throughout the paper, we will study examples whenever they are easily available. Proofs will only be sketched to give an idea of the methods and techniques. We will, however, always include a reference for further reading. After a brief introduction to K3 surfaces, we will study the topics in the above order.

Many of the results and ideas that we explain can be formulated over arbitrary number fields or finite fields. For simplicity we will  mostly restrict to the cases of $\Q$ and $\F_p$.

\section{K3 surfaces}
\label{s:K3}

A.~Weil introduced the notion of K3 surface for any (smooth) surface that carries the differentiable structure of a smooth quartic surface in $\PP^3$. Abstractly a K3 surface is a two-dimensional Calabi-Yau manifold:

\begin{Definition}
A K3 surface is a smooth irreducible surface $X$ with trivial canonical bundle and vanishing first cohomology:
\[
\omega_X\cong\OO_X,\;\;\;h^1(X, \OO_X)=0.
\]
\end{Definition}

The definition allows non-algebraic K3 surfaces. By a theorem of Siu \cite{Siu}, every complex K3 surface is K\"ahler. Any two K3 surfaces are deformation equivalent and thus diffeomorphic. Hence Weil's notion and the above definition coincide. In this survey, we will only consider algebraic K3 surfaces. For details, the reader could consult \cite{BHPV}.

\begin{Example}
A smooth quartic surface in $\PP^3$ is K3. Here we can also allow isolated ADE-singularities. Then the minimal resolution is K3. 
\end{Example}

Quartic surfaces in $\PP^3$ have 19-dimensional moduli. This gives one of countably many components of the moduli space of algebraic K3 surfaces. The component is determined by the polarisation $H^2=4$ where $H$ is an ample line bundle, i.e.~here $H$ corresponds to the hyperplane section.

By Serre duality, a K3 surface $X$ has Euler characteristic $\chi(\OO_X)=2$. We can compute the Euler number with Noether's formula:
\[
e(X) = 12\, \chi(\OO_X) - K_X^2 = 24.
\]
Since $e(X)$ can also be written as alternating sum of Betti numbers, we obtain the Hodge diamond of a K3 surface with entries $h^i(X, \Omega_X^j)$:
$$
\begin{array}{ccccc}
&& 1 &&\\
& 0 && 0 &\\
1 && 20 && 1\\
& 0 && 0 &\\
&& 1 &&
\end{array}
$$
The terminology K3 surface supposedly refers to K\"ahler, Kodaira and Kummer (and to the mountain K2). Kummer surfaces are relevant for many arithmetic aspects of K3 surfaces. Later we will also use elliptic K3 surfaces.

\begin{Example}[Kummer sufaces]
Let $A$ be an abelian surface. Denote the involution by $\iota$. Then the quotient $A/\iota$ has 16 $A_1$ singularities corresponding to the 2-division points on $A$. The minimal resolution is a K3 surface, called Kummer surface Km$(A)$.
\end{Example}

\section{Modularity}
\label{s:mod}

As Calabi-Yau varieties, K3 surfaces are two-dimensional analogues of elliptic curves. In the arithmetic setting, the question of modularity comes to mind. For elliptic curves over $\Q$, modularity is the subject of the Taniyama-Shimura-Weil conjecture as proven by Wiles, Taylor et al.~\cite{Wi}, \cite{TW}, \cite{BCDT}.

\begin{Theorem}[Taniyama-Shimura-Weil conjecture]
Any elliptic curve $E$ over $\Q$ is modular: There is a newform $f$ of weight 2 with Fourier coefficients $a_p$ such that for almost all $p$
\begin{eqnarray}\label{eq:E-mod}
\# E(\F_p) = 1 - a_p + p.
\end{eqnarray}
\end{Theorem}

A newform is a holomorphic function on the upper half plane in $\C$, which satisfies a certain transformation law and is an eigenform of the algebra of Hecke operators. We shall only use that a newform admits a normalised Fourier expansion
\[
\sum_{n\geq 1} a_n q^n,\;\;\;q=e^{2\pi i \tau},\;\; \tau\in\C, \;\mbox{im}(\tau)>0
\]
such that the Fourier coefficients are multiplicative. The complexity of a newform is measured by the level $N\in\N$. Geometrically, the primes dividing the level appear as bad primes for any associated smooth projective variety over $\Q$.

To generalise modularity of elliptic curves over $\Q$, we interpret (\ref{eq:E-mod}) through the Lefschetz fixed point formula. Here we reduce the defining equation of $E$ modulo some prime  $p$. For almost all $p$, this defines a smooth elliptic curve $E_p$. $E_p$ is endowed with the Frobenius morphism Frob$_p$ which raises coordinates to their $p$-th powers. Then the set of $\F_p$-rational points is identified as
\[
E_p(\F_p)=\text{Fix}(\text{Frob}_p).
\]
The induced action Frob$_p^*$ on the cohomology of $E_p$ is subject to the Weil conjectures (proven in generality by Deligne and Dwork). Here one works with $\ell$-adic \'etale cohomology at some prime $\ell\neq p$ which compares nicely between $E$ and $E_p$ (or strictly speaking between the base extensions to an algebraic closure). In the following we will drop the subscripts for simplicity. The Lefschetz fixed point formula yields
\[
\#E(\F_p) = \sum_{i=0}^2 (-1)^i \,\text{trace}(\text{Frob}_p^*;\; H^i(E)).
\]
Then (\ref{eq:E-mod}) is equivalent to $a_p =$ trace$($Frob$_p^*;\; H^1(E))$.

The above ideas generalise directly to any smooth projective variety and its good reductions. For a K3 surface $X$ over $\Q$ and a prime of good reduction $p$, we obtain
\begin{eqnarray}\label{eq:Lef}
\# X(\F_p) = 1 + \text{trace}(\text{Frob}_p^*;\;H^2(X)) + p^2.
\end{eqnarray}
For modularity we have newforms with Fourier coefficients $a_p\in\Z$ in mind. In consequence, modularity requires two-dimensional Galois representations (like those associated to $H^1(E)$). However, $h^2(X)=22$ for a K3 surface. Here we  distinguish that $H^2(X)$ contains both algebraic and transcendental cycles.

\section{Algebraic and transcendental cycles}

The divisors on any smooth projective surface up to algebraic equivalence form the N\'eron-Severi group $\NS(X)$. The rank of $\NS(X)$ is called the Picard number $\rho(X)$. 
Here we can consider $X$ and $\NS(X)$ over any given base field.
Unless specified otherwise, we will be concerned with the geometric N\'eron-Severi group of the base change $\bar X$ to an algebraic closure of the base field -- so that $\NS(X)$ is independent of the chosen model.
On a K3 surface, algebraic and numerical equivalence coincide; hence it suffices to compute intersection numbers. 

$\NS(X)$ is always generated by divisor classes that are defined  over a finite extension of the base field. Hence the absolute Galois group acts by permutations. It follows that all eigenvalues of Frob$_p^*$ are roots of unity. The eigenvalues can be computed explicitly from generators of $\NS(X)$. 

On a complex K3 surface $X$, cup-product endows $H^2(X, \Z)$ with the structure of the unique even unimodular lattice of signature $(3,19)$. Here $\NS(X)$ embeds primitively with signature $(1,\rho(X)-1)$. In particular, $\rho(X)\leq 20$ which also follows from Lefschetz' theorem. We define the transcendental lattice $T(X)$ as 
\[
T(X) = \NS(X)^\bot\subset H^2(X,\Z).
\]
If $X$ is defined over some number field, the lattices of algebraic and transcendental cycles give rise to Galois representations of dimension $\rho(X)$ resp.~$22-\rho(X)$. Modularity requires a two-dimensional Galois representation that does not factor through a finite representation like $\NS(X)$. Hence we need $\rho(X)=20$.

\begin{Example}[Fermat quartic]\label{Ex:F-1}
Consider the Fermat quartic surface
\[
S=\{x_0^4+x_1^4+x_2^4+x_4^4=0\}\subset\PP^3.
\]
$S$ contains 48 lines where $x_0^4+x_i^4=x_j^4+x_k^4=0,\;\{i,j,k\}=\{1,2,3\}$. Their intersection matrix has rank 20. Hence $\rho(S)=20$ if $S$ is considered over $\C$.
\end{Example}

\section{Singular K3 surfaces}
\label{s:sing}

Complex K3 surfaces of Picard number $\rho=20$ are called singular (in the sense of exceptional, but not non-smooth). In many ways, they behave like elliptic curves with complex multiplication (CM), i.e.~elliptic curves $E$ with End$(E)\supsetneq\Z$. 

CM elliptic curves are fully described in terms of class group theory since End$(E)$ is always an order in an imaginary-quadratic field. Analytically, this identification will be exhibited explicitly in (\ref{eq:E}). For CM elliptic curves, analytic and algebraic theory show a particularly nice interplay. Deuring used this to prove that CM elliptic curves are associated to certain Hecke characters. This in fact implies modularity for the 13 CM elliptic curves over $\Q$. The $j$-invariants of CM elliptic curves are often called singular -- thus the terminology for K3 surfaces.

For a singular K3 surface $X$, the relation to class group theory becomes evident when we express $T(X)$ through its intersection form
\begin{eqnarray}\label{eq:Q}
Q(X)=
\begin{pmatrix}
2a & b\\
b & 2c
\end{pmatrix}.
\end{eqnarray}
The quadratic form $Q(X)$ is unique up to conjugation in $SL_2(\Z)$. We denote its discriminant by $d=b^2-4ac<0$.

\begin{Theorem}[Pjatecki\u\i -\v Sapiro - \v Safarevi\v c, Shioda - Inose]\label{Thm:SI}
The map $X \mapsto Q(X)$ is a bijection from isomorphism classes of singular K3 surfaces to even integral positive-definite binary quadratic forms up to conjugation in $SL_2(\Z)$.
\end{Theorem}

The injectivity follows from the Torelli theorem \cite{PS}. To prove the surjectivity, Shioda and Inose \cite{SI} start with an abelian surface $A$ with $\rho(A)=4$ and given quadratic form $Q$ on the transcendental lattice $T(A)$. By earlier work of Shioda-Mitani \cite{SM}, $A$ is isomorphic to the product of two isogenous elliptic curves $E, E'$ with CM in $K=\Q(\sqrt{d})$. As complex tori $E_\tau=\C/(\Z+\tau\Z)$ in terms of the coefficients of $Q$ as in (\ref{eq:Q}), we have
\begin{eqnarray}\label{eq:E}
E=E_\tau,\;\;\tau = \dfrac{-b+\sqrt{d}}{2a},\;\;\;\;\; E'=E_{\tau'},\;\;\tau' = \dfrac{b+\sqrt{d}}2.
\end{eqnarray}
It follows that the Kummer surface of $A$ is a singular K3 surface with intersection form $2Q$. %on $T(\text{Km}(A))$. 
To obtain a K3 surface with the original intersection form $Q$, Shioda-Inose exhibit a particular elliptic fibration on $\mbox{Km}(A)$. Then a suitable quadratic base change of the base curve is again K3 with $\rho=20$ and intersection form $Q$. The corresponding deck transformation is a Nikulin involution, i.e.~an involution with eight isolated fixed points that leaves the holomorphic 2-form invariant. This construction -- Kummer quotient and Nikulin involution yielding the same K3 surface -- is often called Shioda-Inose structure. By work of Morrison \cite{Mo}, it also applies to certain K3 surfaces of Picard number $\rho\geq 17$.
$$
\begin{array}{ccccc}
A &&&& X\\
& \searrow && \swarrow &\\
&& \text{Km}(A) &&
\end{array}
$$
The Shioda-Inose structure implies that any singular K3 surface can be defined over a number field: By class field theory, $E$ and $E'$ can be defined over the ring class field $H(d)$ of discriminant $d$. The Kummer quotient respects the base field, and the elliptic fibrations are defined over some finite extension. An explicit model for $X$ over $H(d)$ was given in \cite{S-fields}. We will see in section \ref{s:NS} how $H(d)$ relates to $\NS(X)$.

\begin{Example}[Fermat quartic cont'd]\label{Ex:F-2}
On the Fermat quartic surface, one can easily single out 20 lines whose intersection form has determinant $-64$. It was a long standing conjecture that these lines generate $\NS(S)$ in characteristic zero. The first complete proof goes back to Inose \cite{I}. He showed that $T(S)$ is four-divisible as an even integral lattice. The only possibility with discriminant dividing $-64$ is
\[
Q(S)=\begin{pmatrix} 8 & 0\\0 & 8\end{pmatrix}.
\]
By (\ref{eq:E}), we also obtain two further models for $S$: through the Shioda-Inose construction for $E_i, E_{4i}$ or as Kummer surface of $E_i\times E_{2i}$. Note that the latter model is defined over $\Q$ as both elliptic curves are, while the former lives only over $\Q(\sqrt{2})$, the real subfield of $H(-64)=\Q(i, \sqrt{2})$.
\end{Example}

Over some extension of the base field (where all the above maps and generators of $\NS(X)$ are defined), the Shioda-Inose structure determines the zeta function of $X$. The essential factor corresponding to $T(X)$ comes from the Hecke character of the elliptic curves.

\section{Modularity of singular K3 surfaces over $\Q$}

For a singular K3 surface $X$ over a number field, the transcendental lattice $T(X)$ gives rise to a system of two-dimensional Galois representations. We have seen that over some extension of the base field, these Galois representations are related to Hecke characters.  If $X$ is defined over $\Q$, modularity was proven by Livn\'e based on a general result about motives with CM \cite{L}.

\begin{Theorem}[Livn\'e]\label{Thm:mod}
Let $X$ be a singular K3 surface over $\Q$ with discriminant $d$. Then $X$ is modular: There is a newform $f$ of weight 3 with CM by $\Q(\sqrt{d})$ and Fourier coefficients $a_p$ such that for almost all $p$
\begin{eqnarray}\label{eq:a_p}
\text{trace}(\text{Frob}_p^*;\;T(X)) = a_p.
\end{eqnarray}
\end{Theorem}

Let $K=\Q(\sqrt{d})$. The CM property is encoded in the Fourier coefficients at the good primes:\begin{eqnarray}\label{eq:f}
a_p = 
\begin{cases}
\pm 2\, (x^2 - d\, y^2), & \text{if $p$ splits in $K$ and } p^2=x^2+d\, y^2\; (x,y\in\Q^*);\\
0, & \text{if $p$ is inert in $K$.}
\end{cases}
\end{eqnarray}

\begin{Example}[Fermat quartic cont'd]\label{F-3}
The associated newform is sensitive to the precise model over $\Q$. For instance, the Fermat quartic $S$ as given in Ex.~\ref{Ex:F-1} is associated to the unique newform with CM in $\Q(\sqrt{-1})$ and Fourier coefficients
\[
a_p=2\,(x^2-4\,y^2)\;\;\; \text{if } p^2=x^2+4\,y^2\; (x,y\in\N).
\]
This is the newform of weight 3 and level $16$ in \cite[Tab.~1]{S-CM}. Other models of $S$ over $\Q$ (as in (\ref{eq:F-4})) lead to twists of the newform, cf.~the discussion below.
\end{Example}

Mazur and van Straten independently asked the converse question whether every such newform of weight 3 is associated to a singular K3 surface over $\Q$.

We are in the opposite situation than for elliptic curves: A classical construction by Eichler and  Shimura associates an elliptic curve over $\Q$ to any newform of weight two with Fourier coefficients $a_p\in\Z$. Meanwhile modularity was open for a long time.

For singular K3 surfaces over $\Q$ we know modularity by Thm.~\ref{Thm:mod}. However, the only known geometric correspondence for newforms of higher weight is through fibre products of universal elliptic curves due to Deligne \cite{Deligne}. Weight 3 leads to elliptic modular surfaces \cite{ShEMS}. Only for a few newforms, modular elliptic surfaces are K3.

In weight 3, the geometric realisation problem is nonetheless approachable because any newform with Fourier coefficients $a_p\in\Z$ has CM by a result of Ribet \cite{R}. Up to twisting, these newforms are in bijective correspondence with imaginary-quadratic fields whose class group is only two-torsion by the classification in \cite{S-CM}. Here twisting refers to modifying the Fourier coefficients by a Dirichlet character $\chi$:
\[
f=\sum_{n\geq 1} a_n q^n\;\; \mapsto \;\; f\otimes\chi = \sum_{n\geq 1} a_n \chi(n) q^n.
\]
The twist is a newform of possibly different level. Twists can be realised geometrically. We have mentioned this for the Fermat surface in Ex.~\ref{F-3}. It is also evident on double covers such as elliptic curves and surfaces: Let $c\in\Q^*$ squarefree and $g\in{\Q}[x_1,\hdots,x_n]$. Then the following double covers are isomorphic over $ \Q(\sqrt{c})$:
\begin{eqnarray}\label{eq:twist}
\{y^2 = g(x_1,\hdots, x_n)\} \mapsto \{ c\,y^2 = g(x_1,\hdots, x_n)\}.
\end{eqnarray}
On the level of modular forms, this corresponds to the twist by the Jacobi symbol $(\frac c\cdot)$.

There are 65 known imaginary-quadratic fields with class group exponent two. By Weinberger \cite{Wb}, there is at most one more such field. Under a natural assumption on Siegel-Landau zeroes of $L$-series of odd real Dirichlet characters (which would follow from the extended Riemann hypothesis), the known list is complete. 

\begin{Theorem}[Elkies, Sch\"utt]\label{Thm:geom}
Every known newform of weight 3 with Fourier coefficients $a_p\in\Z$ is associated to a singular K3 surface over $\Q$.
\end{Theorem}

The proof of Thm.~\ref{Thm:geom} is constructive. For each of the 65 known imaginary-quadratic fields with class group exponent two, it gives a singular K3 surface over $\Q$ \cite{ES-K3}. Each surface admits an elliptic fibration over $\Q$. Hence we can realise all quadratic twists by (\ref{eq:twist}). The next section gives further details of the construction.

\begin{Remark}\label{Rem:ab}
Theorem \ref{Thm:geom} fails for abelian surfaces. Because $H^1(A)$ carries information about the ring class field $H(d)$, we can only realise newforms for fields with class number one or two in abelian surfaces over $\Q$. However, the class number can be as big as 16 for $\Q(\sqrt{-1365})$. 
Since a  K3 surface has trivial $H^1$ by definition, there are only milder obstructions. One of them still involves $H(d)$ and will be given in~Thm.~\ref{Thm:NS} as a consequence of modularity and the Artin-Tate conjecture. Another obstruction involves lattice theory and restricts the genus of $T(X)$ \cite[Thm.~5.2]{S-fields}.
\end{Remark}

\section{Families of K3 surfaces with high Picard number}
\label{s:fam}

To prove Thm.~\ref{Thm:geom}, Elkies and the author considered one-dimensional families of K3 surfaces over $\Q$ with $\rho\geq 19$ \cite{ES-K3}. By moduli theory, such a family has infinitely many specialisations with $\rho=20$ over $\bar\Q$. We aim at specialisations over $\Q$.

The Lefschetz fixed point formula (\ref{eq:Lef}) provides a good test whether a K3 surface is modular: If $\rho(X)\geq 19$, then we obtain 19 eigenvalues of Frob$_p^*$ on $H^2(X)$ from the Galois operation on $\NS(X)$. By the Weil conjectures, the non-real eigenvalues come in complex conjugate pairs. Hence there is one further eigenvalue $\pm p$. If $\rho=20$, then the remaining two eigenvalues give the trace of Frob$_p^*$ on the Galois representation associated to $T(X)$. By counting points, we can check whether this trace agrees with the coefficient $a_p$ of any given newform of weight 3.

Applied to several primes, this test can rule out or provide evidence for $\rho=20$. But then we have to prove that $\rho=20$ at some specialisation. We sketch two ways to approach this problem.

On the one hand, the parametrising curve is always a modular curve or a Shimura curve. Hence, if we can determine the curve and its CM points, we are done. However, the K3 surfaces that we are interested in become more and more complicated. We will see in Thm.~\ref{Thm:NS} that this can be measured by the discriminant of $X$ or equivalently by the Galois action on $\NS(X)$. But this means that the Shimura curves become actually so complicated that we barely have any knowledge about them at all. In fact, families of K3 surfaces with $\rho\geq 19$ provide a new tool to find detailed information about otherwise inaccessible Shimura curves (cf.~\cite{E-Shimura}).

\begin{Example}[Dwork pencil]
The Fermat quartic can be deformed into several one-dimensional familes with $\rho\geq 19$.
Consider the famous Dwork pencil
\[
S_\lambda=\{x_0^4+x_1^4+x_2^4+x_4^4=\lambda\, x_0x_1x_2x_3\}\subset\PP^3,\;\;\;\;\; \lambda\neq 0.
\]
One way to determine the parametrising curve is related to mirror symmetry: The quotient $Y_\lambda$ of $S_\lambda$ by a $(\Z/4\Z)^2$ subgroup of $\mbox{Aut}(S_\lambda)$ is a family of K3 surfaces with $\rho=19$ and  discriminant $4$. By \cite{Dol}, the parametrising Shimura curve is $X_0(2)^+$. The parametrising curve of $S_\lambda$ is a four-fold cover of $X_0(2)^+$.
A new algebraic approach to this problem using Shioda-Inose structures is pursued in \cite{ES}.
\end{Example}

On the other hand, we can look for an additional divisor on some specialisation which would imply $\rho=20$. This problem seems untractable for the Dwork pencil. It becomes feasible when turning to elliptic surfaces with section. This also has the side-effect that we can twist as in (\ref{eq:twist}).

The N\'eron-Severi group of an elliptic surface is generated by fiber components and sections due to the Shioda-Tate formula. Hence to increase $\rho$, the singular fibers could degenerate, but in general there has to be an independent section. To find a specialisation with an independent section $P$, we have to solve a system of at least seven non-linear equations.

A newform that we want to realise geometrically fixes the discriminant of $\NS(X)$ up to square by Theorem \ref{Thm:mod}. The theory of Mordell-Weil lattices after Shioda \cite{ShMW} translates $\mbox{disc}(\NS(X))$ into conditions on the intersection numbers of $P$ with fibre components and the other sections. The fibre component restrictions cut down the number of equations we have to solve. Sometimes, this suffices to determine a solution directly. In other cases, we exhibit an extensive seach over some finite field $\F_p$ to find a special K3 surface with an independent section mod $p$. Then we employ a $p$-adic Newton iteration to increase the $p$-adic accuracy. Finally we compute a lift to $\Q$ with the Euclidean algorithm  and verify that $\rho=20$ for this specialisation.

\section{Picard numbers of surfaces}
\label{s:Pic}

The Picard number of a projective surface is in general hard  to determine. A method by Shioda applies to surfaces with large group actions, such as Fermat surfaces and their quotients \cite{Sh-Pic}. In the previous section, we used reduction mod $p$ and the Lefschetz trace formula to check for modularity. By Thm.~\ref{Thm:mod}, this is equivalent to $\rho=20$ (cf.~Ex.~\ref{Ex:19}).  These ideas will be put in a systematic context in the sequel. Then we will give some applications.

Let $X$ be a smooth projective surface over $\Q$. Then $X$ has good reduction $X_p$ at almost all primes $p$. From now on, we have to distinguish between the geometric N\'eron-Severi group $\NS(\bar X)$ of the base extension $\bar X$ to an algebraic closure and the sublattice $\NS(X)$ generated by divisors over the given base field. The reduction morphism induces embeddings of lattices
\begin{eqnarray}\label{eq:NS}
\NS(X) \hookrightarrow \NS(X_p),\;\;\;\;\;\; \NS(\bar X)\hookrightarrow \NS(\bar X_p).
\end{eqnarray}
In particular, the Picard number cannot decrease upon good reduction.

We have already noted that Frob$_p^*$ operates as identity on $\NS(X_p)$. On $\NS(\bar X_p)$, there are roots of unity involved which we shall denote by $\zeta$. 
On the algebraic subspace of $H^2(X)$, the eigenvalues of Frob$_p^*$ are exactly these roots of unity multiplied by $p$.

The eigenvalues of Frob$_p^*$ on $H^2(X)$ are encoded in the characteristic polynomial
\begin{eqnarray}\label{eq:char}
R_p(T) = \det(T-\text{Frob}_p^*; H^2(X)).
\end{eqnarray}
By Newton's identities for symmetric polynomials, one obtains all coefficients of $R_p(T)$ from the traces of the powers of Frob$_p^*$. Hence it suffices to count points and apply the Lefschetz fixed point formula (\ref{eq:Lef}) for sufficiently many extensions $\F_q$ where we replace $p$ by $q=p^s$. 
The computations are simplified by Poincar\'e duality and the fact that $R_p(T)$ is monic. Hence, for a K3 surface, one at worst has to consider $\F_{p^{11}}$.

%By the Weil conjectures, the eigenvalues on $H^2$ have absolute value $p$ and come in complex conjugate pairs. Thus one only has to determine half as many coefficients of $R_p(T)$. For a K3 surface, one at worst has to consider $\F_{p^{11}}$. 

The eigenvalues of Frob$_p^*$ on $H^2(X)$ give upper bounds for the Picard numbers:
\begin{eqnarray}\label{eq:Tate}
\rho(X_p)\leq \text{ord}_{T=p}(R_p(T)),\;\;\; \rho(\bar X_p)\leq \sum_\zeta \text{ord}_{T=\zeta p}(R_p(T)).
\end{eqnarray}
Sometimes these bounds suffice to compute the Picard number of a surface: when we compute some divisors and their rank equals the bound at some good prime $p$ obtained from (\ref{eq:Tate}). 

However, this approach alone cannot work for all surfaces in characteristic zero. The reason lies in a parity condition due to the Weil conjecture: Since all eigenvalues of Frob$_p^*$ on $H^2(X)$ other than $\pm p$ come in complex conjugate pairs, the following differences are even:
\begin{eqnarray}\label{eq:b2}
b_2(X) - \text{ord}_{T=\pm p}(R_p(T)) ,\;\; b_2(X) - \sum_\zeta \text{ord}_{T=\zeta p}(R_p(T)).
\end{eqnarray}

\begin{Conjecture}[Tate {\cite{Tate-C}}]\label{Conj:T}
The inequalities in (\ref{eq:Tate}) are in fact equalities.
\end{Conjecture}

Tate proved the conjecture for abelian surfaces and products of curves. By work of Artin and Swinnerton-Dyer on principal homogeneous spaces \cite{ASD}, the Tate conjecture also holds for elliptic K3 surfaces with section.

The Tate conjecture predicts the parity of $\rho(\bar X_p)$ by (\ref{eq:b2}). Hence, if a K3 surface over $\Q$ has odd Picard number, then $\rho$ has to increase upon reduction. Therefore we cannot use the above method directly to compute $\rho(\bar X)$. We sketch an idea due to van Luijk to circumvent this parity problem in the next section.

\section{Van Luijk's method}
\label{s:van}

Van Luijk \cite{van} pioneered a method to prove odd Picard number on a K3 surface $X$ over $\Q$. The setup is as follows: Assume we have a lower bound $\rho(\bar X)\geq r$ with odd $r$, say from explicit divisors on $X$. Find a prime $p$ such that the eigenvalues of Frob$_p^*$ on $H^2(X)$  imply $\rho(\bar X_p)\leq r+1$ by (\ref{eq:Tate}). Then van Luijk's method gives a criterion to show that the lower bound is attained.

Assume that $\rho(\bar X)=r+1$. Then the embeddings in (\ref{eq:NS}) are isometries up to finite index. Hence the discriminants of the N\'eron-Severi groups in characteristic zero and $p$ agree up to the square of the index. In this sense, the discriminants have to be compatible for all reductions at good primes $p$ with $\rho(\bar X_p)=r+1$.

\begin{Criterion}[van Luijk]\label{Thm:van}
In the above setup, assume that there are good primes $p_1\neq p_2$ such that $\rho(\bar X_{p_i})= r+1$. %by (\ref{eq:Tate}). 
Let 
\[
D=\dfrac{\mbox{disc}(\NS(\bar X_{p_1}))}{\mbox{disc}(\NS(\bar X_{p_2}))}.
\]
If $D$ is not a square in $\Q^*$, then $\rho(\bar X)\leq r$.
\end{Criterion}

\begin{Example}\label{Ex:19}
In section \ref{s:fam} we implicitly used this technique to prove $\rho=19$ for K3 surfaces in a one-dimensional family over $\Q$: Otherwise $\rho(X)=20$ and $X$ would be modular by Thm.~\ref{Thm:mod}. 
The point counting test checks whether (\ref{eq:a_p}) holds at any given $p$ for a newform $f$ of weight 3. Here $f$ has CM by $K=\Q(\sqrt{d})$, encoded in the Fourier coefficients. 

At the inert primes in $K$, this coefficient is zero. This results in further eigenvalues $p, -p$ of Frob$_p^*$ on $H^2(X)$. Subject to the Tate conjecture, the Picard number increases upon reduction, so we cannot use the above criterion. At the split primes in $K$, however, the Picard number stays constant upon reduction by (\ref{eq:Tate}). Since $K$ is determined by $d$ up to square, the test amounts to checking Criterion~\ref{Thm:van}.
\end{Example}

Criterion~\ref{Thm:van} requires the discriminant of the N\'eron-Severi group at two good primes up to square. In practice, this does not mean that one has to compute explicit divisors and their intersection numbers. Kloosterman noticed that instead one can work with the Artin-Tate conjecture which we formulate here for K3 surfaces over $\F_q$. It involves the characteristic polynomial $R_q(T)$ of Frob$_q^*$ on $H^2(X)$ as in (\ref{eq:char}).

\begin{Conjecture}[Artin-Tate {\cite{Artin-Tate}}]\label{Conj:AT}
Let $X$ be a K3 surface over $\F_q$. Let $\alpha(X_q)=b_2(X)-\rho(X_q)-1$. Then
\begin{eqnarray}\label{eq:AT}
\dfrac{R_q(T)}{(T-q)^{\rho(X_q)}}\Big |_{T=q} =  q^{\alpha(X_q)}\; |Br(X_q)|\cdot |\text{discr}(\NS(X_q))|.
\end{eqnarray}
\end{Conjecture}

By \cite{Milne} (cf.~\cite[p.~25]{Milne-A} for characteristic two), the Artin-Tate conjecture is equivalent to the Tate conjecture. For us, it is crucial that the order of the Brauer group $Br(X_q)$ is always a square by \cite{Brauer}. Hence we can compute the discriminant of the N\'eron-Severi group up to square by $R_q(T)$. By the Lefschetz fixed point formula (\ref{eq:Lef}), this amounts to point counting over sufficiently many extensions of $\F_q$. Note that if $X$ is defined over $\F_p$, then we obtain $R_q(T)$ from $R_p(T)$ by Newton's identities for any $q=p^s$. Hence we only have to compute $R_p(T)$, working over the corresponding extensions of $\F_p$.

We now return to the setup of a K3 surface $X$ over $\Q$ with $\rho(\bar X)\geq r$ for some odd $r$. Based on Crit.~\ref{Thm:van}, we give an algorithm to prove that $\rho(\bar X)=r$. The third step which replaces the actual computation of $\mbox{disc}(NS(X_q))$ is due to Kloosterman \cite{Kl}.

\begin{enumerate}[1.]
\item Find a good prime $p$ and the characteristic polynomial $R_{p}(T)$ such that $\rho(\bar X_p)\leq r+1$ by (\ref{eq:Tate}).
\item Choose $q=p^s$ such that $R_{q}(T)$ has the root $T=q$ with multiplicity $r+1$. 
\item Define $\tilde D_q = \dfrac 1q \dfrac{R_q(T)}{(T-q)^{r+1}}\Big |_{T=q}$. 
\item Repeat the above procedure for another good prime $p$.
\item If the $\tilde D_q$ are not multiples by a square factor, then $\rho(\bar X)\leq r$ by Crit.~\ref{Thm:van}.
\end{enumerate}

\begin{Remark}
In the above algorithm, we do not have to assume the Tate conjecture: In order to establish a contradiction, we only make the assumption that $\rho(\bar X)=r+1$. By the choice of $q$, the Tate conjecture follows automatically from the specialisation embedding (\ref{eq:NS}). Hence the Artin-Tate conjecture gives the discriminant up to squares. We compute it in the third step and compare with the discriminant for other suitable choices of $q$. Thus we can derive a contradiction without assuming the Tate conjecture.
\end{Remark}

\section{Applications of van Luijk's method}
\label{s:app}

We have already seen one application of van Luijk's method for families of K3 surfaces of Picard number $\rho\geq 19$. The method was independently based on the modularity of singular K3 surfaces over $\Q$ (cf.~Ex.~\ref{Ex:19}). In this section, we will sketch two other applications, due to van Luijk and Kloosterman.

\subsection{K3 surfaces with Picard number one}

By moduli theory, a general complex algebraic K3 surface has Picard number one. 
Terasoma showed that there is a smooth quartic K3 surface over $\Q$ with $\rho=1$ \cite{Terasoma}.
This is a non-trivial fact since there are only countably many K3 surfaces over $\Q$. 
However, not even over $\C$, there was a single explicit K3 surface with $\rho=1$ known until van Luijk formulated Crit.~\ref{Thm:van} in \cite{van}.

Here the main challenge is computational: The determination of $R_p(T)$ can require point counting over fields $\F_q$ with $q=p^s, s=1,\hdots 11$. Hence van Luijk looked for a K3 surface $X$ over $\Q$ with good reduction at the two smallest primes, 2 and 3. Then he computed $R_p(T)$.

At each prime, the second requirement is that $R_p(T)$ has only two roots $\zeta p$. Once this is achieved, the discriminants of the N\'eron-Severi groups are used to establish the claim $\rho(X)=1$. In fact, van Luijk obtained a much stronger result:

\begin{Theorem}[van Luijk]\label{Thm:1}
The smooth quartic K3 surfaces over $\Q$ with $\rho=1$ are dense in the moduli space of K3 surface with a polarisation of degree 4.
\end{Theorem}

The proof relies on the above K3 surface $X$ which van Luijk had found as a smooth quartic, say
\[
X=\{f=0\}\subset\PP^3,\;\;\; f\in{\Q}[x_0,x_1,x_2,x_3] \text{ homogenous of degree } 4.
\]
Now consider a family of K3 surfaces parametrised in terms of another homogenous polynomial $h\in{\Q}[x_0,x_1,x_2,x_3]$ of degree $4$:
\begin{eqnarray}\label{eq:X_h}
\frak X_h = \{f=6\,h\}\subset\PP^3.
\end{eqnarray}
Whenever $h$ has coefficients in $\Q$ that are integral $2$-adically and $3$-adically, then $\frak X_h$ reduces to the same surfaces $X_2$ mod $2$ and $X_3$ mod $3$. Hence $\rho=1$ for all these surfaces by the above argument. The surfaces are easily seen to lie dense in the given component of the moduli space of K3 surfaces.

%\subsection{A 17-dimensional family of K3 surfaces over $\Z$ with $\rho=3$}

%For any 17-dimensional family of K3 surfaces, moduli theory gives the upper bound $\rho\leq 3$ for a general member. Baragar and van Luijk construct a particular family over $\Z$ such that $\rho\geq 3$.

%The remarkable fact about this family is that for all parameters in $\Z$, the member has Picard number $\rho=3$. To achieve this, the family is chosen in such a way that all these surface reduce to the same surfaces $X_2$ and $X_3$ mod $2$ resp.~$3$.
%Then one applies the above methods to these single surfaces to deduce that $\rho(X_2), \rho(X_3)\leq 4$. Since $D_2/D_3$ is not a square in $\Q^*$, the claim $\rho(X)\leq 3$ follows for any lift to $\Z$. By construction, this applies to the whole family over $\Z$. (Note that parameter choices in $\Q$ can very well lead to increasing Picard number.)

\subsection{An elliptic K3 surface with Mordell-Weil rank 15}
\label{s:15}

A complex elliptic K3 surface can have Mordell-Weil rank from 0 to 18, since $\rho\leq 20$. Cox proved that all these ranks occur \cite{Cox}. This result is a consequence of the surjectivity of the period map, but Cox did not give any explicit examples. As a complement, Kuwata derived explicit complex elliptic K3 surfaces over $\Q$ with any Mordell-Weil rank from 0 to 18 except for 15 \cite{Kuwata}. Kloosterman filled the gap with help of the algorithm we described in the previous section \cite{Kl}.

Through various base changes involving rational elliptic surfaces, Kloosterman derived an elliptic K3 surface with Mordell-Weil rank at least $15$. Then he proved equality by the above algorithm. Since $\rho\geq 17$, he knew a factor of $R_p(T)$ of degree 17. Hence the determination of $R_p(T)$ only required point counting over $\F_p$ and $\F_{p^2}$. Therefore the technique was also feasible at the larger primes $p=17, 19$ which Kloosterman used.

\section{Class group action on singular K3 surfaces}
\label{s:NS}

The Artin-Tate conjecture has many other applications to K3 surfaces over finite fields or number fields. Here we sketch one of them -- which again does not require to assume the validity of the conjecture.

To prove Thm.~\ref{Thm:geom}, we constructed singular K3 surfaces over $\Q$ where the associated imaginary-quadratic field $K=\Q(\sqrt{d})$ has class number 4, 8 or even 16. At first, this might come as a surprise since the equivalent statement for abelian surfaces does not hold because of the relation to the class field $H(d)$ (cf.~Rem.~\ref{Rem:ab}). On the other hand, we know from class field theory that both surfaces can be defined over $H(d)$ (cf.~sect.~\ref{s:sing}). In this section, we will see that independent of the field of definition, every singular K3 surface carries the class group structure: through the Galois action on the N\'eron-Severi group.

\begin{Theorem}[Elkies, Sch\"utt]\label{Thm:NS}
Let $X$ be a singular K3 surface. Let $L$ be a number field such that $\NS(\bar X)$ is generated by divisors over $L$. Denote the discriminant of $X$ by $d$. Let $H(d)$ be the ring class field for $d$. Then $H(d)\subset L(\sqrt{d})$.
\end{Theorem}

Here we shall only sketch the case where $X/\Q$ is a singular K3 surface with $\NS(\bar X)$ generated by divisors over $\Q$. This case is originally due to Elkies \cite{Elkies}. Thm.~\ref{Thm:NS} can be rephrased as follows:

\begin{Theorem}[Elkies]\label{Thm:Q}
Let $X/\Q$ be a singular K3 surface. Assume that $\NS(\bar X)$ is generated by divisors over $\Q$. Then $X$ has discriminant $d$ of class number one.
\end{Theorem}

We shall sketch an alternative proof given by the author in \cite{S-NS} which uses modularity and the Artin-Tate conjecture. The main idea of the proof is as follows:

Consider the good primes $p$ that split in $K=\Q(\sqrt{d})$. By Thm.~\ref{Thm:mod},
\[
R_p(T) = T^2- a_p\,T + p^2.
\]
Since $p\nmid a_p$, the reduction $X_p$ has $\rho(X_p)=\rho(\bar X_p)=20$ by (\ref{eq:Tate}). Hence the Artin-Tate conjecture holds for $X_p$. From (\ref{eq:AT}) we obtain a relation between $a_p$ and $d$ up to squares. First we ignore the squares and consider only $K$. Using the explicit description of $a_p$ in (\ref{eq:f}), one can show that $p$ splits into principal ideals in $K$. Since by assumption this holds for almost all split $p$, $K$ has class number one.

Then we take the precise shape of $d$ into account. This is made possible by the fact, that the embeddings (\ref{eq:NS}) are almost always primitive. With elementary class group theory using representations of numbers by quadratic forms, one can prove that $d$ also has class number one.

\begin{Remark}
In Thm.~\ref{Thm:Q}, one can also show that $T(X)$ is primitive. Otherwise the singular K3 surface would be Kummer or admit an isotrivial elliptic fibration with $j=0$, but Mordell-Weil rank two.

Conversely, for each $d<0$ with class number one, there is a singular K3 surface $X$ with $\NS(\bar X)$ generated by divisors over $\Q$ and discriminant $d$ (cf.~\cite[\S 10]{S-NS}).
\end{Remark}

We will study an application of Thm.~\ref{Thm:Q} to elliptic K3 surfaces in the next section. Here we only note that Thm.~\ref{Thm:NS} gives a direct proof of \v Safarevi\v c' finiteness theorem for singular K3 surfaces which generalises the theory for CM elliptic curves:

\begin{Theorem}[\v Safarevi\v c {\cite{Shafa}}]
Let $n\in\N$. Up to $\C$-isomorphism, there are only finitely many singular K3 surfaces over number fields of degree at most $n$.
\end{Theorem}

It is an open problem to determine all singular K3 surfaces over a fixed number field, say over $\Q$. It follows from work of Shioda \cite{Sh-base}, that the absolute value of the discriminant can be at least as big as $36\cdot 427$.

\section{Ranks of elliptic curves}
\label{s:ranks}

In section \ref{s:15}, we have cited that every Mordell-Weil rank from 0 to 18 is possible on complex elliptic K3 surfaces. The question of the rank over $\Q$ is much more delicate. Shioda asked in \cite{Shioda-20} whether Mordell-Weil rank 18 over $\Q$ is possible. Elkies gave a negative answer in \cite{Elkies}, \cite{Elkies-rank}, based on Thm.~\ref{Thm:Q}.

Since the N\'eron-Severi group of an elliptic surface is generated by horizontal and vertical divisors, Mordell-Weil rank 18 implies that an elliptic K3 surface is singular. Moreover all fibres have to be irreducible. In consequence, the Mordell-Weil lattice is even and integral, but has no roots (i.e.~elements with minimal norm $2$). This already is a special property.

If the Mordell-Weil rank over $\Q$ is 18, then $\NS(\bar X)$ is generated by divisors over $\Q$. Hence Thm.~\ref{Thm:Q} bounds the discriminant $|d|\leq 163$. Because of the absence of roots, such a lattice would break the density records for sphere packings in $\R^{18}$. By gluing up to a  Niemeier lattice, Elkies is able to establish a contradiction.

On the other hand, Elkies found an elliptic K3 surface with Mordell-Weil rank 17 over $\Q$ (and necessarily $\rho=19$) \cite{Elkies-rank}. All intermediate ranks are also attained. After base changing from this elliptic K3 surface, a certain specialisation yields an elliptic curve with rank (at least) 28 over $\Q$. This extends the previous record by 4.

\section{Rational points on K3 surfaces}

Given a variety $X$ over a field $k$, it is natural to ask for the $k$-rational points $X(k)$. This set can be empty, finite, infinite or even dense in $X$. For instance, the model of the Fermat quartic in Ex.~\ref{Ex:F-1} has no rational points over any totally real field. In this section we will review the situation for K3 surfaces over number fields.

The most general case would be K3 surfaces of Picard number one. In 2002, Poonen and Swinnerton-Dyer asked whether there is a K3 surface with $\rho=1$ over a number field with infinitely many rational points. Van Luijk gave a striking affirmative answer in \cite{van}:

\begin{Theorem}[van Luijk]
In the moduli space of K3 surface with a polarisation of degree 4,
the K3 surfaces over $\Q$ with $\rho=1$ and infinitely many rational points are dense.
\end{Theorem}

Van Luijk's basic idea is to find a dense set of K3 surfaces within the family $\frak X_h$ in (\ref{eq:X_h}) which contain an elliptic curve of positive rank over $\Q$. By construction, the rational points thus obtained are not dense on the respective K3 surface. The question of density of rational points thus remains open. It has been answered for some other K3 surfaces with $\rho>1$. We mention one example due to Harris and Tschinkel {\cite{Harris-Tschinkel}}:

\begin{Theorem}[Harris, Tschinkel]
Let $S$ be a smooth quartic in $\PP^3$ over some number field $k$. Assume that $S$ contains a line $\ell$ over $k$ which does not meet more than five lines on $S$. Then $S(k)$ is dense in $S$.
\end{Theorem}

The proof uses the fact that $S$ admits an elliptic fibration. In Thm.~\ref{Thm:BT} we will see more general consequences of this property in the context of potential density. 

Recently Logan, van Luijk and McKinnon announced a similar result for twists of the Fermat quartic. They consider those models
\begin{eqnarray}\label{eq:F-4}
S'=\{a\,x_0^4+b\,x_1^4+c\,x_2^4+d\,x_3^4=0\}\subset\PP^3
\end{eqnarray}
with $a,b,c,d\in\Q^*$ such that $abcd$ is a square. They assume that there is a rational point on $S$ outside the 48 lines and the coordinate planes. Then they conclude that the $\Q$-rational points are dense on $S$.

Here we touch on two crucial problems:
\begin{enumerate}[1.]
\item Is there a rational point on a given projective variety?
\item Does the existence of a rational point (plus possibly some conditions) imply that there are infinitely many rational points? Will the rational points be dense?
\end{enumerate}

The first question motivated the Hasse principle: If $X$ has a $\Q$-rational point, then it has $\Q_v$-rational points at every place $v$. Hence the set of adelic points $X(\A)$ is nonempty. The converse implication, known as the Hasse principle
\[
X(\A)\neq\emptyset \;\, \stackrel{?}{\Longrightarrow} \;\, X(\Q)\neq \emptyset,
\]
holds for conics. It is computationally very useful since the existence of local points can be checked algorithmically. However, the Hasse principle is not true in general. Selmer showed that it fails for the genus one curve
\[
C=\{3\,x^3+4\,y^3+5\,z^3=0\}\subset\PP^2.
\]
In 1970, Manin \cite{Manin} discovered that the failure of the Hasse principle can often be explained through the Brauer group $\Br(X)$. He defined a subset $X(\A)^{\Br}\subset X(\A)$ that contains $X(\Q)$, but can be empty even if $X(\A)$ is not. This case is referred to as Brauer-Manin obstruction. 

It is conjectured that the Brauer-Manin obstruction is the only obstruction to the Hasse principle on rational varieties. More generally, however, Skorobogatov \cite{Sko} constructed a bielliptic surface where the Brauer-Manin obstruction does not suffice to explain the failure of the Hasse principle. Recent developments indicate that even natural cohomological refinements of the Brauer-Manin obstruction may not be sufficient in general \cite{Poonen}. 

For K3 surfaces, it is still an open problem whether the Brauer-Manin obstruction is the only obstruction to the Hasse principle. For the models of the Fermat surface over $\Q$ in (\ref{eq:F-4}), Swinnerton-Dyer \cite{SD} proved this under the following assumptions: Schinzel's hypothesis holds, the Tate-\v Safarevi\v c group of elliptic curves is finite plus conditions on the coefficients $a, b, c, d$. Here the conditions on the coefficients guarantee the existence of an elliptic fibration over $\Q$. Then the elaborate techniques from \cite{CTSSD} apply (cf.~\cite[pp.~585, 625/626]{CTSSD} for an account of arithmetic implications).

The second problem is currently being investigated by van Luijk. Numerically he found evidence for the rate of growth of the number of points with bounded height $h$, thus indicating an affirmative answer. The conjectural formula involves the Picard number, resembling the situation for del Pezzo surfaces as predicted by Manin's conjecture. In some instances, the formula requires to leave out a finite number of curves on the surface. On a Zariski-open subset $U$, the conjectural formula in case $\rho(X)=1$
reads 
\[
\# \{P\in U(\Q); h(P)\leq B\} \sim c\cdot\log B, \;\;\;\;\;\; c \text{ constant}.
\]

Rational points and their density rely heavily on the model of the variety and the chosen base field -- like for the Fermat surface in Ex.~\ref{Ex:F-1} and its twists in (\ref{eq:F-4}). This suggests that the notion of density would be too restrictive for most general statements. This is the reason to introduce potential density.

\section{Potential density}

Let $X$ be a variety over a number field or the function field of a curve, denoted by $k$. We say that the rational points  are  \emph{potentially dense} on $X$ if there is some finite extension $k'/k$ such that $X(k')$ is dense. The main expectation and motivation for this notion is that potential density should be a geometric property, depending only on the canonical class $K_X$.

For curves, this concept is known thanks to Faltings' theorem \cite{F}: For any curve of genus greater than one over any number field, the rational points are finite. Conversely, potential density holds for curves of genus zero and one.

The same is true for surfaces with $K_X$ negative: Over a finite extension, they are all rational. On the other end, the Lang-Bombieri conjecture rules out potential density for projective varieties of general type.

Among surfaces with $K_X\equiv 0$, potential density has been proved for abelian and Enriques surfaces. A result by Bogomolov and Tschinkel covers a great range of K3 surfaces \cite{BT}.

\begin{Theorem}[Bogomolov, Tschinkel]\label{Thm:BT}
Let $X$ be a K3 surface over a number field. Assume that $X$ has an elliptic fibration or infinite automorphism group. Then the rational points  are potentially dense on $X$.
\end{Theorem}

Here either assumption implies $\rho\geq 2$. Other than this, the assumptions are not too restrictive. For instance, any K3 surface with $\rho=2$, but without $(-2)$ curves satisfies the conditions of the theorem. Similarly, any K3 surface with $\rho\geq 5$ admits an elliptic fibration and therefore shares the property of potential density.

It is again the case $\rho=1$ where we lack any examples over number fields -- be it with or without potential density of rational points. The only known result concerns function fields of complex curves:

\begin{Theorem}[Hassett, Tschinkel {\cite{HT}}]
Let $C$ be a complex curve. There are non-trivial K3 surfaces over $\C(C)$ with $\rho=1$ and Zariski-dense rational points.
\end{Theorem}

The proof relies on the uncountability of $\C$. It is unclear how the techniques could be adapted for function fields $\bar\Q(C)$. Despite the recent progress, the question of potential density  is still wide open for K3 surface.

\vspace{0.5cm}

\textbf{Acknowledgement:} 
This survey owes to the contributions of many mathematicians from whom I learned through lectures, discussions and collaborations.
My thanks go especially to N.~D.~Elkies, B.~van Geemen, K.~Hulek, R.~Kloosterman, R.~van Luijk, T.~Shioda, and to the referee. 
Funding from DFG 
under research grants Schu 2266/2-1 and Schu 2266/2-2 is gratefully acknowledged.

\vspace{0.5cm}

Matthias Sch\"utt,
Department of Mathematical Sciences,
University of Copenhagen,
Universitetspark 5,
2100 Copenhagen,
Denmark,
{\tt mschuett@math.ku.dk}

\end{document}